\documentclass{article}
\usepackage{amsrefs}
\usepackage{amssymb}
\usepackage{amsmath}
\usepackage{amsthm}
\usepackage{graphicx}

\DeclareMathOperator{\sepgon}{sepgon}
\DeclareMathOperator{\Proj}{Proj} 
\DeclareMathOperator{\Ext}{Ext} \numberwithin{equation}{section}
\newtheorem{theorem}{Theorem}
\newtheorem{definition}[theorem]{Definition}

\newtheorem{lemma}[theorem]{Lemma}
\newtheorem{definition/construction}[theorem]{Definition/Construction}

\theoremstyle{remark}
\newtheorem*{remark}{Remark}
\newtheorem*{subject}{2000 Mathematics Subject Classification}
\newtheorem*{keywords}{Keywords}

\author{Marc Coppens\footnote{Katholieke Hogeschool Kempen, Departement IBW,
Kleinhoefstraat 4, B-2440 Geel, Belgium; K.U.Leuven, Departement of Mathematics,
Celestijnenlaan 200B, B-3001 Leuven, Belgium; email: marc.coppens@khk.be}}
\title{The separating gonality of a separating real curve}
\date{}

\begin{document}
\maketitle \noindent

\begin{abstract}

Let $X$ be a smooth real curve of genus $g$ such that the real locus
has $s$ connected components. We say $X$ is separating if the
complement of the real locus is disconnected. In case there exists a
morphism $f$ from $X$ to $\mathbb{P}^1$ such that the inverse image
of the real locus of $\mathbb{P}^1$ is equal to the real locus of
$X$ then $X$ is separating and such morphism is called separating.
The separating gonality of a separating real curve $X$ is the
minimal degree of a separating morphism from $X$ to $\mathbb{P}^1$.
It is proved by Gabard that this separating gonality is between $s$
and $(g+s+1)/2$. In this paper we prove that all values between $s$
and $(g+s+1)/2$ do occur.

\end{abstract}

\begin{subject}
14H05; 14H51; 14P99
\end{subject}

\begin{keywords}
real curve, gonality, rational function, linear system
\end{keywords}

\section{Introduction}

Let $X$ be a smooth real curve of genus $g$. We assume $X$ is
complete and geometrically irreducible, hence its set
$X(\mathbb{C})$ of complex points is in a natural way a compact
Riemann surface of genus $g$. Let $X(\mathbb{R})$ be the set of real
points, it is well known that $X(\mathbb{C})\setminus X(\mathbb{R})$
is either connected or it has two connected components. In case it
has two connected components then we say $X$ is a \emph{separating
real curve}.

It is well known that each component of $X(\mathbb{R})$ is a smooth
analytic real manifold of dimension 1 homeomorphic to a circle. Let
$s=s(X)$ be the number of connected components of $X(\mathbb{R})$.
In case $X$ is separating one knows $1\leq s\leq g+1$ and $s \equiv
g+1 \pmod{2}$. For each integer $s$ satisfying those conditions
there exists a real separating curve $X$ of genus $g$ with $s(X)=s$.

Let $f:X\rightarrow \mathbb{P}^1$ be a morphism such that
$f^{-1}(\mathbb{P}^1(\mathbb{R}))=X(\mathbb{R})$. Hence each real
fiber of $f$ only consists of real points, such morphism is called
\emph{totally real} or \emph{separating}. We use the terminology of
separating: since $\mathbb{P}^1(\mathbb{C})\setminus
\mathbb{P}^1(\mathbb{R})$ is not connected the existence of such
morphism $f$ implies $X$ is separating. Conversely, it is proved in
\cite{ref5} that for each separating real curve $X$ of genus $g$
there exists a separating morphism of degree at most $g+1$. In
\cite{ref6} one proves a stronger result: each separating real curve
$X$ of genus $g$ with $s=s(X)$ has a separating morphism of degree
at most $(g+s+1)/2$. At the end of the introduction of that paper,
the author mentions that sharpness of that upper bound is an open
problem. Related to this problem and copying the terminology of
gonality of a complex curve, we introduce the following definition.

\begin{definition}
Let $X$ be a separating real curve of genus $g$. The
\emph{separating gonality} of $X$ is the minimal degree $k$ such
that there exists a separating morphism $f:X \rightarrow
\mathbb{P}^1$ of degree $k$. This separating gonality is denoted by
$\sepgon (X)$.
\end{definition}

The result in \cite{ref6} can be written as $\sepgon (X)\leq
(g+s+1)/2$ with $s=s(X)$. Moreover, there is the trivial lower bound
$\sepgon (X)\geq s$. In this paper we prove the following theorem,
giving amongst others an affirmative answer to the question on the
sharpness on the upper bound for $\sepgon$.

\begin{theorem}
Let $g$, $s$ and $k$ be integers satisfying $g\geq 2$, $1\leq s\leq
g+1$ with $s \equiv g+1 \pmod{2}$ and $\max \{ 2,s \}\leq k\leq
(g+s+1)/2$. Then there exists a separating real curve $X$ of genus
$g$ with $s(X)=s$ and $\sepgon (X)=k$.
\end{theorem}

In case $s=g+1$ then the lower and upper bound on $\sepgon$ agree.
Such curves are called M-curves and separating morphisms of degree
$g+1$ on M-curves are studied in \cite{ref8}. As a matter of fact,
if $C_1, \cdots, C_{g+1}$ are the connected components of the real
locus of an M-curve $X$ and $P_i \in C_i$  for $1\leq i\leq g+1$
then $\{ P_1, \cdots, P_{g+1} \}$ is the fiber of a separating
morphism $f:X\rightarrow \mathbb{P}^1$, uniquely determined by those
points up to an automorphism of $\mathbb{P}^1$. However the other
extreme case $s=1$ is the base case for the proof of the theorem. In
this case the upper bound on $\sepgon$ is equal to the Meis bound
for the gonality of compact Riemann surfaces (\cite{ref9}) and this
is nowadays contained in Bril-Noether Theory (see e.g.
\cite{ref10}*{Chapter V}). Using results concerning pencils on
Riemann surfaces we obtain a proof of the theorem in case $s=1$
(this argument resembles those used in \cite{ref7}*{Section 3}).
Then we take a suited separating real curve $X'$ with $s(X')=1$ and
on $X'_{\mathbb{C}}$ we identify closed points associated to some
chosen non-real points on $X'$. This gives rise to a singular real
curve $X_0$ on which each identification gives rise to an isolated
real node. We consider a suited smoothing $X_t$ of $X_0$ defined
over $\mathbb{R}$ such that each isolated real node becomes a
component of $X_t(\mathbb{R})$ and such that $X_t$ has a suited
separating morphism to $\mathbb{P}^1$. Such smoothing $X_t$ is a
separating real curve with $s(X_t)=s$ and we prove that if $X_t$ is
close enough to $X_0$ then it has the desired separating gonality.

\section{Notations and preliminaries}

Let $\mathbb{R}$ be the field of real numbers, let $X$ be an
indeterminate over $\mathbb{R}$ and let $\mathbb{C}$ be the field
$\mathbb{R}[X]/(X^2+1)$, called the field of complex numbers. The
complex number defined by $X$ is denoted by $i$; in this way we have
the field extension $\mathbb{R}\subset \mathbb{C}$. For a complex
number $z=x+iy$ ($x$,$y\in \mathbb{R}$), its complex conjugate is
denoted by $\overline{z}=x-iy$.

A real scheme is a reduced, separated scheme of finite type over
$\mathbb{R}$. In case $V$ is a real scheme then $V_{\mathbb{C}}$ is
the scheme obtained from $V$ by the base change defined by
$\mathbb{R}\subset \mathbb{C}$. We say $V$ is smooth, irreducible,
... in case $V_{\mathbb{C}}$ is smooth, irreducible, ... . We write
$V(\mathbb{C})$ to denote the set of closed points of
$V_{\mathbb{C}}$. On $V_{\mathbb{C}}$, hence also on
$V(\mathbb{C})$, we have the $\mathbb{R}$-morphism induced by
complex conjugation on $\mathbb{C}$. It is called complex
conjugation on $V_{\mathbb{C}}$ and on $V(\mathbb{C})$. In case
$P\in V(\mathbb{C})$ then $\overline{P}$ is the point obtained by
complex conjugation. There are two types of closed points on $V$:
those corresponding to $P\in V(\mathbb{C})$ with $P=\overline{P}$
and those corresponding to $\{ P, \overline{P} \}\subset
V(\mathbb{C})$ with $P\neq \overline{P}$. In case $P=\overline{P}$
then $P$ is called a \emph{real point} on $V$. The set of real
points on $V$ is denoted by $V(\mathbb{R})$, it is a subset of
$V(\mathbb{C})$. In case $P\neq \overline{P}$ then the associated
point on $V$ is called a non-real point on $V$ and it is denoted by
$P+\overline{P}$. In case $V$ and $W$ are two real schemes then a
morphism $f:V\rightarrow W$ is defined over $\mathbb{R}$. We write
$f_{\mathbb{C}}:V_{\mathbb{C}}\rightarrow W_{\mathbb{C}}$ to denote
the morphism defined by the field extension $\mathbb{R}\subset
\mathbb{C}$. It induces maps $f(\mathbb{C}):V(\mathbb{C})\rightarrow
W(\mathbb{C})$ and $f(\mathbb{R}):V(\mathbb{R})\rightarrow
W(\mathbb{R})$. We write $\mathbb{P}^N$ do denote the scheme $\Proj
(\mathbb{R}[X_0, \cdots, X_N])$ (see \cite{ref11}). It is called the
real projective space of dimension $N$. In accordance to the general
notations then $\mathbb{P}^N_{\mathbb{C}}$ is the complex projective
space of dimension $N$ and we have an inclusion
$\mathbb{P}^N(\mathbb{R}) \subset \mathbb{P}^N(\mathbb{C})$.

For some basic terminology used now we refer to \cite{ref11}. A real
scheme $X$ is called a real curve if $X_{\mathbb{C}}$ is a
connected, reduced curve (hence it has dimension 1 but it can be
singular); $X$ is called complete, irreducible, smooth, stable, ...
in case $X_{\mathbb{C}}$ is complete, irreducible, smooth, stable
... . Now assume $X$ is a real complete curve.  An
$\mathcal{O}_X$-Module $L$ that is locally isomorphic to
$\mathcal{O}_X$ is called an invertible sheaf on $X$. Using the
field extension $\mathbb{R} \subset \mathbb{C}$ we obtain an
invertible sheaf $L_{\mathbb{C}}$ on $X_{\mathbb{C}}$. The space of
global sections $\Gamma (X,L)$ (resp. $\Gamma
(X_{\mathbb{C}},L_{\mathbb{C}})$) is a finite dimensional
$\mathbb{R}$-vectorspace (resp. $\mathbb{C}$-vectorspace). Complex
conjugation on $X$ induces a complex conjugation on $\Gamma
(X_{\mathbb{C}},L_{\mathbb{C}})$ having invariant space $\Gamma
(X,L)$. In this way we consider $\Gamma
(X_{\mathbb{C}},L_{\mathbb{C}})$ as the complexification of $\Gamma
(X,L)$ and we denote $h^0(L)$ for both $\dim _{\mathbb{R}}(\Gamma
(X,L))$ and $\dim _{\mathbb{C}}(\Gamma
(X_{\mathbb{C}},L_{\mathbb{C}}))$. A Cartier divisor $D$ on $X$
defines an invertible sheaf $\mathcal{O}_X(D)$ and it is well-known
that $\mathbb{P} (\Gamma (X, \mathcal{O}_X (D)))$ (resp. $\mathbb{P}
(\Gamma (X_{\mathbb{C}}, \mathcal{O}_{X_{\mathbb{C}}} (D)))$)
parameterizes the space of effective divisors on $X$ (resp.
$X_{\mathbb{C}}$) linearly equivalent to $D$. This is called the
complete linear system defined by $D$, it is denoted by
$|D|=\mathbb{P} (\Gamma (X, \mathcal{O}_X (D)))$ and
$|D|_{\mathbb{C}}=\mathbb{P} (\Gamma (X_{\mathbb{C}},
\mathcal{O}_{X_{\mathbb{C}}} (D)))$. A linear subspace $V$ of
$\Gamma (X, \mathcal{O}_X (D))$ defines a linear subsystem
$\mathbb{P}(V)$ of $|D|$. In case $\deg (D)=d$ and $r+1=\dim
_{\mathbb{R}} (V)$ then we say $\mathbb{P}(V)$ is a $g^r_d$ on $X$.
On $X_{\mathbb{C}}$ we have
$\mathbb{P}(V)_{\mathbb{C}}=\mathbb{P}(V_{\mathbb{C}})$ (here
$V_{\mathbb{C}}\subset \Gamma (X_{\mathbb{C}},
\mathcal{O}_{X_{\mathbb{C}}} (D))$ is the complexification of $V$)
and $g^r_{d,\mathbb{C}}$. Of course using subvectorspaces of the
global space of an invertible sheaf $L$ on $X_{\mathbb{C}}$ we
obtain linear systems on $X_{\mathbb{C}}$ (but not necessarily
defined by linear systems on $X$).

The moduli functor of stable curves of genus $g$ is not
representable, hence there is no fine module scheme with a universal
family. Instead we make use of so-called suited families of stable
curves.

\begin{definition}

Let $X$ be a complete stable real curve of genus $g$. A \emph{suited
family of stable curves of for $X$} is a projective morphism $\pi :
\mathcal{C}\rightarrow S$ defined over $\mathbb{R}$ such that

\begin{enumerate}
\item
$S_{\mathbb{C}}$ is smooth, irreducible and quasi-projective.

\item
Each fiber of $\pi_{\mathbb{C}}$ is a stable curve of genus $g$.

\item
For each $s\in S(\mathbb{C})$ the Kodaira-Spencer map
$T_s(S_{\mathbb{C}}) \rightarrow \Ext ^1(\Omega
_{\pi_{\mathbb{C}}^{-1}(s)},\mathcal{O}_{\pi_{\mathbb{C}}
^{-1}(s)})$ is surjective (here $\Omega _{\pi_{\mathbb{C}}
^{-1}(s)}$ is the sheaf of K\"ahler differentials).

\item
There exists $s_0 \in S(\mathbb{R})$ such that $X$ is isomorphic to
$\pi^{-1}(s_0)$.
\end{enumerate}

\end{definition}

From the third condition it follows that for a suited family
$\pi:\mathcal{C} \rightarrow S$ of stable curves of genus $g$ one
has $\dim (S)=3g-3$ in case $g\geq 2$.

\begin{lemma}
Let $X$ be a complete stable curve of genus $g\geq 2$. There exists
a suited family of stable curves for $X$.
\end{lemma}

\begin{proof}
We are going to use some facts mentioned in \cite{ref13}*{Section
1}). The moduli functor of tricanonically embedded stable curves of
genus $g$ is represented by a quasi-projective scheme $\mathcal{H}$
defined over $\mathbb{R}$ and there exists $x\in
\mathcal{H}(\mathbb{R})$ corresponding to a tricanonical embedding
of $X$. This scheme is smooth and moreover the natural
transformation to the moduli functor of stable curves obtained by
forgetting the embedding of the fibers is formally smooth. Hence the
induced Kodaira-Spencer map $T_x(\mathcal{H})\rightarrow
\Ext^1(\Omega _X,\mathcal{O}_X)$ (the tangent space of the moduli
functor of stable curves at $X$) is surjective. Taking an embedding
of $\mathcal{H}$ in some projective space and intersecting with
suited hyperplanes through $x$ one obtains a suited family of stable
curves for $X$.
\end{proof}

As alread mentioned in the introduction, we are going to use
singular real curves obtained by making real isolated singular ponts
out of some fixed non-real points. First we prove a lemma on the
canonical embedding of those singular curves that we need in the
proof. A lot of geometry of a canonically embedded curve is
determined by the Riemann-Roch Theorem, so first we recall this
theorem.

Let $X$ be a real nodal curve of arithmetic genus $g$. It has a
dualizing sheaf $\omega _X$ (defined over $\mathbb{R}$). This is an
invertible sheaf on $X$ and it satisfies the Riemann-Roch Theorem:
for an invertible sheaf $L$ on $X$ one has

\begin{equation*}
h^0(L)=\deg (L)-g+1+h^0(\omega _X\otimes L^{-1})
\end{equation*}

\noindent (see e.g. \cite{ref1}*{Chapter 4}). In particular $\omega
_X$ is the unique invertible sheaf of degree $2g-2$ whose space of
global sections has dimension $g$. Now let $1\leq m\leq g-2$ be an
integer and let $X'$ be a smooth real curve of genus $g'=g-m$ and
for $1\leq i\leq m$ let $P_i + \overline{P_i}$ be different non-real
points on $X'$ such that $\dim |P_i + \overline{P_i}| =0$ (such
non-real points exist because $X'$ has at most one $g^1_2$ and this
cannot contain all non-real points). Let $X$ be the singular nodal
real curve obtained from $X'$ by identifying $P_i$ with
$\overline{P_i}$ for $1\leq i\leq m$. This identification gives rise
to an isolated real singular point $S_i$ on $X$. The arithmetic
genus of $X$ is $g$.

\begin{lemma}\label{lemmaveryample}
$\omega _X$ is very ample.
\end{lemma}

\begin{proof}
For $1\leq j\leq m$ let $X_j$ be the curve obtained from $X'$ by
identifying $P_i$ with $\overline{P_i}$ for $1\leq i\leq j$ and
again write $S_i$ for the resulting nodes on $X_j$. One has $X_m=X$
and we write $X_0$ instead of $X'$. So for each $1\leq j\leq m$ the
curve $X_j$ is obtained from the curve $X_{j-1}$ by identifying
$P_j$ and $\overline{P_j}$ considered as points on $X_{j-1}$. First
we prove $\omega _{X_1}$ is very ample on $X_1$ and then, taking
$2\leq j\leq m$ and assuming $\omega _{X _{j-1}}$ is very ample on
$X_{j-1}$ we prove $\omega _{X_j}$ is very ample on $X_j$.

From the Riemann-Roch Theorem on $X_{0}$ it follows $h^0(\omega
_{X_0}(P_1+\overline{P_1}))=g'+1$ and for all $Q\in X_0(\mathbb{C})$
one has $h^0(\omega _{X_0,\mathbb{C}}(P_1+\overline{P_1}-Q))=g'$.
This proves the complete linear system $|\omega
_{X_0}(P_1+\overline{P_1})|_{\mathbb{C}}$ on $X_{0,\mathbb{C}}$ is
base point free, hence it defines a morphism $i_1 : X_0 \rightarrow
\mathbb{P}^{g'}$ having image $i_1(X_0)$. Since $h^0(\omega
_{X_0,})=g'=h^0(\omega _{X_0}(P_1+\overline{P_1}))-1$ it follows
$i_1(\mathbb{C})(P_1)=i_1(\mathbb{C})(\overline{P_1})=S'_1$ is a
real singular point on $i_1(X_0)$. Since $g'\neq 0$ for $Q\in
X_0(\mathbb{C})$ one has $h^0(\mathcal{O}_{X_0,\mathbb{C}}(Q))=1$,
hence $h^0(\omega _{X_0,\mathbb{C}}(-Q))=g-1=h^0(\omega
_{X_0}(P_1+\overline{P_1}))-2$. This proves
$i_1(\mathbb{C})^{-1}(S'_1)=\{ P_1, \overline{P_1} \}$ and
$i_1(\mathbb{C})$ is a local immersion at both $P_1$ and
$\overline{P_1}$. Projection from $\mathbb{P}^{g'}$ to
$\mathbb{P}^{g'-1}$ with center $S'_1$ induces the canonical map
$i_0:X_0\rightarrow \mathbb{P}^{g'-1}$. Since
$P_1+\overline{P_1}\notin g^1_2$ (if $X_0$ would have a $g^1_2$) it
follows $i_0(\mathbb{C})(P_1)\neq i_0(\mathbb{C})(\overline{P_1})$,
hence $S'_1$ is an ordinary node on $i_1(X_0)$. In case
$i_1(X_0)_{\mathbb{C}}$ would have another singular point $S'$ then
there exists an effective divisor $F$ of degree 2 on
$X_0(\mathbb{C})$ having disjoint support with $P_1+\overline{P_1}$
such that $h^0(\omega _{X_0,\mathbb{C}}(P_1+\overline{P_1}-F))=g'$.
Because of the Riemann-Roch Theorem it would imply $\omega
_{X_0,\mathbb{C}}(P_1+\overline{P_1}-F)\cong \omega
_{X_0,\mathbb{C}}$ hence $F$ and $P_1+\overline{P_1}$ are linearly
equivalent. This again contradicts $P_1+\overline{P_1}\notin g^1_2$.
So we obtain $i_1(X_0)=X_1$ and $X_1\subset \mathbb{P}^{g'}$ is
embedded by an invertible sheaf on $X_1$ of degree $2g'=2p_a(X_1)-2$
since $p_a(X_1)=g'+1$. It follows $X_1\subset \mathbb{P}^{g'}$ is
canonical embedded, hence $\omega _{X_1}$ is very ample.

Now, let $2\leq j\leq m$ and assume $\omega _{X_{j-1}}$ is very
ample. Let $i_{j-1}:X_{j-1}\rightarrow \mathbb{P}^{g'+j-2}$ be a
canonical embedding. Consider $\omega
_{X_{j-1}}(P_j+\overline{P_j})$ on $X_{j-1}$. From the Riemann-Roch
Theorem on $X_{j-1}$ it follows the complete linear system $|\omega
_{X_{j-1}}(P_j+\overline{P_j})|$ is base point free, hence it
defines a morphism $i_j:X_{j-1}\rightarrow \mathbb{P}^{g'+j-1}$. One
finds again $i_j(P_j)=i_j(\overline{P_j})=S'_j$ is a real singular
point on $i_j(X_{j-1})$ and projection of $\mathbb{P}^{g'+j-1}$ to
$\mathbb{P}^{g'+j-2}$ with center $S'_j$ induces the canonical
embedding $i_{j-1}$ of $X_{j-1}$. This implies $S'_j$ is an ordinary
node and $i_j(\mathbb{C})$ defines an isomorphism between
$X_{j-1}\setminus \{ P_j, \overline{P_j} \}$ and
$i_j(X_{j-1})\setminus \{ S'_j \}$. This proves $X_j\cong
i_j(X_{j-1})$ and as before we conclude $X_j \subset
\mathbb{P}^{g'+j-1}$ is canonically embedded.
\end{proof}

\begin{remark}
This lemma and its proof corresponds to the well-known description
of the dualizing sheaf  on a nodal curve as it is described in e.g.
\cite{ref1}*{Chapter IV.9}. According to that description a section
of the dualizing sheaf of $X'$ corresponds to a rational
differential form on the normalization $X$ having poles of order at
most one on the inverse images of the nodes such that the sum of
their residues at those 2 points in 0 and being regular outside the
nodes. So the sections belong to $\Gamma (X, \omega _{X}(\sum
_{i=1}^m (P_i+\overline{P_i})))$ and if such section is 0 at some
$P_i$ (or $\overline{P_i}$) (meaning it is regular as differential
form at $P_i$ (or $\overline{P_i}$)) then because of the restiction
on the residus, it need to be regular as differential forms on both
$P_i$ and $\overline{P_i}$, hence the section vanishes at both $P_i$
and $\overline{P_i}$.

\end{remark}

\section{Proof of the theorem}

First we are going to prove the theorem in case $s=1$. Since $X$ is
separated it follows $g(X)$ is even, say $2h$. In this case $\sepgon
(X)=s=1$ is impossible hence $2\leq \sepgon(X)\leq (g+2)/2=h+1$.
This upper bound is equal to the upper bound on the gonality of
complex curves of genus $2h$. For $2\leq k\leq h+1$ it is proved in
\cite{ref7}*{Theorem 2} that there exists a real curve $X$ of
topological type $(2h,1,0)$ having a morphism $f:X\rightarrow
\mathbb{P}^1$ of degree $k$ such that $f(\mathbb{R}):X(\mathbb{R})
\rightarrow \mathbb{P}^1(\mathbb{R})$ is an unramified covering of
degree $k$. For $k\leq h$ one has $2k-2h-2<0$ (hence the
Brill-Noether Number $\rho ^1_k(2h)$ (see \cite{ref10}*{Chapter IV})
is negative) and therefore it follows from \cite{ref7}*{Theorem 9}
that a general real curve $X$ of topological type $(2h,1,0)$ does
not have a separating morphism of degree $k$. This proves the
upperbound on $\sepgon$ obtained in \cite{ref6} is sharp in case
$s=1$. In \cite{ref6}*{Remark 11} it is mentioned that there exist
real curves of topological type $(2h,1,0)$ having a separating
morphism of degree $k$ and not having a separating morphism of
degree less than $k-1$. This finishes in case $s=1$ the proof of our
theorem. Since this case is crucial for the proof of our theorem for
all other cases we give some more details concerning those claims in
\cite{ref7}*{Remark 11}.

\begin{lemma}
Let $h,k$ be integers with $h\geq 1$ and $2\leq k\leq h$. There
exists a real curve $X$ of topological type $(2h,1,0)$ such that $X$
has a separating morphism $f:X\rightarrow \mathbb{P}^1$ of type $k$
and $X_{\mathbb{C}}$ has no linear system $g^1_l$ with $l<k$.
\end{lemma}

\begin{proof}
From \cite{ref7}*{Theorem 2}, it follows there exists a real curve
$X'$ of topological type $(2h,1,0)$ having a separating morphism
$f':X'\rightarrow \mathbb{P}^1$. Let $\pi : \mathcal{C}\rightarrow
S$ be a suited family of curves of genus $g$ for $X'$. Choose $s'\in
S(\mathbb{R})$ with $\pi ^{-1}(s')=X'$. As in \cite{ref7} we
consider the $S$-scheme $\pi _k:H_k(\pi ) \rightarrow S$
parameterizing morphisms of degree $k$ from fibers of $\pi$ to
$\mathbb{P}^1$. It is a smooth quasi-projective variety of dimension
$2g+2k-2$ defined over $\mathbb{R}$. In particular $H_k(\pi
)(\mathbb{R})$ parameterizes real morphisms of degree $k$ from real
fibers of $\pi$ to $\mathbb{P}^1$ and this is a real analytic
manifold (not connected) of real dimension $2g+2k-2$. In $\pi
_k^{-1}(s')$ there is a point $[f']$ on $H_k(\pi )(\mathbb{R})$
corresponding to $f'$. It belongs to a connected component $H$ of
$H_k(\pi )(\mathbb{R})$ and since the topological degree (see
\cite{ref7}) of real morphisms is a discrete continuous invariant,
it is constant on $H$. In our situation this means each point of $H$
represents a separating morphism of degree $k$ on a real curve of
topological type $(2h,1,0)$ (also this topological type cannot
change in a connected family of real curves). The dimension of
$H_k(\pi)$ is in accordance with the Hurwitz formula for the number
of ramification points of a covering
$f_{\mathbb{C}}:X_{\mathbb{C}}\rightarrow \mathbb{P}^1$ of degree
$k$ with $g(X_{\mathbb{C}})=g$ (see e.g. \cite{ref11}*{Chapter IV,
Corollary 2.4}). Such covering is called simple in case all
ramification points have index two and no two of them map to the
same point of $\mathbb{P}^1_{\mathbb{C}}$. A moduli count using
ramification shows that the subspace $H^{ns}_k(\pi )$ parameterizing
non-simple coverings has dimension less than $2g+2k-2$. Since it is
invariant under complex conjugation, it is defined over $\mathbb{R}$
and therefore $\dim H^{ns}_k(\pi ))(\mathbb{R})<2g+2k-2$ hence $H
\nsubseteq H^{ns}_k(\pi )(\mathbb{R})$. Let $l<k$ and consider
$H_{k,l}=H_k(\pi) \times _S H_l(\pi)$ and let $H'_{k,l}(\pi)$ be the
image of $H_{k,l}$ on $H_k(\pi)$. This is a constructible subset of
$H_k(\pi)$ and $H'_{k,l}(\pi)$ is invariant under complex
conjugation, so it is defined over $\mathbb{R}$. Assume
$([g],[h])\in (H_{k,l})_{\mathbb{C}}$ and $g$ is simple. This
defines a morphism $(g,h):X \rightarrow
\mathbb{P}^1_{\mathbb{C}}\times \mathbb{P}^1_{\mathbb{C}}$ (here $X$
is the complex fiber of $\pi$ associated to $[g]$). Assume the image
is not birational equivalent to $X$, then it defines non-trivial
morphisms $g':X \rightarrow X'$ and $g'':X' \rightarrow
\mathbb{P}^1_{\mathbb{C}}$ with $g=g''\circ g'$. Since $g$ is
simple, $g''$ can not have remification, a contradiction to the
Hurwitz formula. So the image of $(g,h)$ is birationally equivalent
to $X$. Then \cite{ref6}*{Proposition 2.4} implies $\dim
_{([g],[h])}(H_{k,l})=g+2k+2l-7$. Since $l\leq h=g/2$ it follows
$\dim _{([g],[h])}(H_{k,l})<2g+2h-2$. This proves
$\dim(H'_{k,l}(\pi))<2g+2k-2$ and so $H\nsubseteq
H'_{k,l}(\pi)(\mathbb{R})$. Take $[f]\in H$ with $[f] \notin
H'_{k,l}(\pi)(\mathbb{R})$ for all $l<k$ and $s=\pi _k([f])\in
S(\mathbb{R})$ then $\pi ^{-1}(s)$ is a real curve of topological
type $(2h,1,0)$ and $f:\pi ^{-1}(s) \rightarrow \mathbb{P}^1$ is a
separating morphism of degree $k$ and $\pi ^{-1}(s)_{\mathbb{C}}$
has no morphism $g:\pi ^{-1}(s)_{\mathbb{C}} \rightarrow
\mathbb{P}^1_{\mathbb{C}}$ of degree $l<k$.
\end{proof}

Now assume $s>1$. For $s\leq g+1$ and $s \equiv g+1 \pmod{2}$ it is
proved in \cite{ref7}*{Theorem 2} that there exists a real curve $X$
of topological type $(g,s,0)$ having a separating morphism $f:X
\rightarrow \mathbb{P}^1$ of degree $s$. Of course, for such curve
we have $\sepgon (X)=s$, hence we can restrict to the case $s+1 \leq
k\leq (g+s+1)/2)$ (in particular also $s\leq g-1$). Let $2h=g-s+1$
and $k'=k-s+1$, hence $2h\geq 2$ and $2\leq k'\leq (g-s+1)/2
+1=h+1$. Take a real curve $X'$ of topological type $(2h,1,0)$ such
that $\sepgon (X')=k'$ and $X'_{\mathbb{C}}$ has no morphism to
$\mathbb{P}^1_{\mathbb{C}}$ of degree $l<k'$. Let $f':X'\rightarrow
\mathbb{P}^1$ be a separating morphism of degree $k'$. Choose
general points $P_2, \cdots, P_s$ on $X'(\mathbb{C})$. Take $s-1$
copies of $\mathbb{P}^1$ denoted by $(\mathbb{P}^1)_2, \cdots,
(\mathbb{P}^1)_s$ and let $\Gamma _0$ be the real curve obtained
from $X' \cup (\mathbb{P}^1)_2 \cup \cdots \cup (\mathbb{P}^1)_s$ by
identifying $P_i$ to $f'(P_i)$ and $\overline{P_i}$ to
$f'(\overline{P_i})=\overline{f'(P_i)}$ for $2\leq i\leq s$. This
singular real curve $\Gamma _0$ has arithmetic genus $g$ and the
morphism $f'$ together with the identities between
$(\mathbb{P}^1)_i$ and $\mathbb{P}^1$ for $2\leq i\leq s$ give rise
to a morphism $f_0:\Gamma _0 \rightarrow \mathbb{P}^1$ defined over
$\mathbb{R}$ of degree $k$. Moreover for $P\in
\mathbb{P}^1(\mathbb{R})$ the fiber $f_0^{-1}(P)$ only has real
points, hence $f_0$ is a separating morphism for this singular real
curve. The associated stable curve $X_0$ is the curve obtained from
$X'$ by identifying $P_i$ and $\overline{P_i}$ for $2\leq i\leq s$,
hence it has a real isolated node $S_i$ for $2\leq i\leq s$ as its
only singularities.

Locally at the image of $P_i$ on $\Gamma _0$ (still denoted by
$P_i$) the curve $\Gamma _0(\mathbb{C})$ has a description $z^2
-w^2$ with $P_1=(0,0)$ and with $f_0(z,w)=z$. Let $U$ be a small
neighborhood of $(0,0)$ in $\mathbb{C}^2$ and $V\subset U$ a much
smaller one. Let $U_0=\Gamma _0(\mathbb{C})\cap U$ and $V_0=\Gamma
_0(\mathbb{C})\cap V$. For $P\neq P_i$ on $U_0$ we can use $z$ to
define a holomorphic coordinate at $P$. Consider a local deformation
$z^2 -w^2=t$ with $t\in \mathbb{C}$ and $|t|$ small and let $U_t$
(resp. $V_t$) be the intersection with $U$ (resp. $V$). We use a
gluing of $U_t$ and $\Gamma _0(\mathbb{C}) \setminus V_0$ as
follows. For $z_0\in \mathbb{C}^*$ let $ _{z_0}\sqrt{z}$ be the
locally defined holomorphic square root function such that $
_{z_0}\sqrt{z_0^2}=z_0$. A point $Q\in U_0 \setminus V_0$ has
coordinates $(z,z)$ or $(z,-z)$. We identify $(z,z)$ with $(z,
_z\sqrt{z^2-t})$ and $(z,-z)$ with $(z,- _z\sqrt{z^2-t})$. (One
should adopt the description of $U_0$ and $V_0$ to this
identification). Moreover the projection $(z,w)\rightarrow z$ on
$U_t$ and $f_0(\mathbb{C})$ restricted to $\Gamma
_0(\mathbb{C}\setminus V_0 )$ glue together to a holomorphic mapping
to $\mathbb{P}^1(\mathbb{C})$. At $\overline{P_i}$ (the image of
$\overline{P_i}$ on $\Gamma _0$) one has a complex conjugated
description $\overline{z}^2-\overline{w}^2=0$ with
$\overline{P_i}=(0,0)$ and
$f_0(\overline{z},\overline{w})=\overline{z}$. Using similar
neighborhoods denoted by $\overline{U}$ and $\overline{V}$ one can
consider a similar local smoothing
$\overline{z}^2-\overline{w}^2=t'$. In case $t'=\overline{t}$ there
is a complex conjugation between $U_t$ and
$\overline{U}_{\overline{t}}$ and it is compatible with the complex
conjugation on $\Gamma _0 \setminus (V_0 \cup \overline{V}_0)$.
Parameterizing this deformation by $(x,y)\in \mathbb{C}^2$ with
$x=t+t'$ and $y=(1/i)(t-t')$ we obtain the deformation has a global
antiholomorphic involution if $(x,y)\in \mathbb{R}^2$. We write
$\Gamma _{(x,y)}$ for the deformation defined by $(x,y)$, it has a
real structure for $(x,y)\in \mathbb{R}^2$. Moreover for $(x,y)\in
\mathbb{R}^2$ the morphism $f_0 : \Gamma _0 \rightarrow
\mathbb{P}^1$ deforms to a morphism $f_{(x,y)} : \Gamma _{(x,y)}
\rightarrow \mathbb{P}^1$. Doing the deformation for $2\leq i\leq s$
we obtain the existence of $\epsilon
>0$ such that there is a deformation $\pi : \mathcal{X} \rightarrow
T= \{ z\in \mathbb{C}: |z|< \epsilon \}$ and a $T$-morphism $F :
\mathcal{X} \rightarrow T \times \mathbb{P}^1$ such that $\pi
^{-1}(0)=\Gamma _0$, $F| _{\Gamma _0} : \Gamma _0 \rightarrow
\mathbb{P}^1$ is $f_0$ and for $t \in T(\mathbb{R})$ the curve $\pi
^{-1}(t)=X_t$ is a real curve of topological type $(g,s,0)$ and $F|
_{X_t}:X_t \rightarrow \mathbb{P}^1$ is a separating morphism of
degree $k$. As a matter of fact, for $P\in \mathbb{P}^1(\mathbb{R})$
the fiber $f_t^{-1}(P)$ has degree $k'$ on the component $C_{1,t}$
of $X_t(\mathbb{R})$ degenerating to $X'(\mathbb{R})$ and it has
degree 1 on each component $C_{j,t}$ of $X_t(\mathbb{R})$
degenerating to the component $(\mathbb{P}^1)_j(\mathbb{R})$ of
$\Gamma _0(\mathbb{R})$. It follows $\sepgon (X_t)\leq k$ and we are
going to prove that $\sepgon (X_t)=k$ for $t\in \mathbb{R}\setminus
\{ 0 \}$ it $|t|$ is sufficiently small.

Let $X_0$ be the stable model of $\Gamma _0$ and let $\pi :
\mathcal{C}\rightarrow S$ be a suited family for $X_0$. let $s_0\in
S(\mathbb{R})$ with $\pi ^{-1}(s_0)=X_0$. Let $\omega _{\pi}$ be the
relative dualizing sheaf of $\pi$. The restriction of $\omega
_{\pi}$ to $\pi ^{-1}(s_0)$ is the dualizing sheaf $\omega_{X_0}$,
it is very ample because of Lemma \ref{lemmaveryample}. Hence
shrinking $S$ we can assume $\omega _{\pi ^{-1}(s)}$ is very ample
for all $s\in S(\mathbb{C})$.  Let $\mathcal{E}=\pi _{*}(\omega
_{\pi})$ a locally free sheaf of rank $g$ on $S$ defined over
$\mathbb{R}$. Shrinking $S$ we can assume $\mathcal{E}$ is free and
taking a base on it defined over $\mathbb{R}$ we obtain an
$S$-morphism $\mathcal{C}\rightarrow S\times \mathbb{P}^{g-1}$ and
the image is a real family of canonically embedded curves. We use
this family and we keep writing $\pi : \mathcal{C} \rightarrow S$
(so now this is a family of canonically embedded curves). From the
previous arguments it follows there exists a real curve $T\subset
S(\mathbb{R})$ containing $s_0$ such that for $t\in T \setminus
\{s_0 \}$ the curve $\pi ^{-1}(t)=X_t$ is a curve of topological
type $(g,s,0)$ and it has a separating morphism $f_t : X_t
\rightarrow \mathbb{P}^1$ of degree $k$. To make notations easier we
write $0$ to denote $s_0$ in the following arguments. We write
$C_{1,t}$ to denote the components of $X_t(\mathbb{R})$ deforming to
the smooth component $C_{1,0}$ of $X_0(\mathbb{R})$ and $C_{i,t}$ to
denote the component of $X_t(\mathbb{R})$ deforming to $S_i$ for
$2\leq i\leq s$.

Assume for all $t\in T\setminus \{ 0 \}$ there exists a separating
morphism $f_t : X_t \rightarrow \mathbb{P}^1$ of some degree $l<k$.
Shrinking $T$ if necessary there exists some $l<k$ such that for
each $t\in T \setminus \{ 0 \}$ there exists a separating morphism
$f_t : X_t \rightarrow \mathbb{P}^1$ of degree $l$. Let $D_t$ be a
fiber of $f_t$, then $L_t=\mathcal{O}(D_t)$ is an invertible sheaf
of degree $l$ satisfying $h^0(L_t)\geq 2$. Since $l\leq g-1$ it
follows from the Riemann-Roch Theorem that $h^0(\omega
_{X_t}(-D_t))\geq 1$. Hence there is a hyperplane section $H$ of
$X_t$ (an effective canonical divisor) containing $D_t$. Let
$E_t=H-D_t$. For each effective divisor $D$ linearly equivalent to
$D_t$ there is a hyperplane section of $X_t$ equal to $D+E_t$. From
$h^0(\mathcal{O} _{X_t}(D_t))\geq 2$, it follows there is a
one-dimensional family of hyperplanes in $\mathbb{P}^{g-1}$
containing $E_t$. The intersection of such hyperplanes is a linear
subspace $\Lambda _t$ of dimension at most $g-3$ containing $E_t$.
(As a matter of fact the Riemann-Roch Theorem implies $\dim (\Lambda
_t)=g-1-h^0(\mathcal{O}(D_t))$.) Those linear subspaces define a
family of linear subspaces over $T \setminus \{ 0 \}$ and this
family has a limit $\Lambda _0$ above 0 (such limit exists because
the Grassmannian is complete). Since $\Lambda _t \cap X_t$ is a
scheme of length $2g-2-l$ for $t\neq 0$, in the limit $\Lambda _0
\cap X_0$ contains a scheme of length $2g-2-l$ (here we use the fact
that the relative Hilbert scheme is proper over the base). Let
$C_{1,0}$ be the smooth component of $X_0(\mathbb{R})$ and let $P\in
C_{1,0}$. Let $P_t$ be a family of points over $T$ such that $P_t\in
C_{1,t}$ for $t\in T$ with $P_0=P$. For $t\neq 0$ the fiber $D_t$ of
$f_t$ containing $P_t$ also contains some point $P_{t,j}$ of
$C_{j,t}$ for $2\leq j\leq s$. Let $H_t$ be the hyperplane
containing $\Lambda _t$ such that $H_t\cap X_t=D_t+E_t$ (this is an
intersection as scheme, hence possibly with multiple points). The
limit of $H_t$ is a hyperplane $H_0$ and the limit of $D_t+E_t$ is
equal to the scheme $H_0\cap X_0$. This subscheme contains $P\in
C_{t,0}$. Moreover the limit of $E_t$ is contained in $\Lambda
_0\cap X_0$ and the limit of $D_t$ contains the limit of the points
$P_{t,j}$ for $2\leq j\leq s$, hence it contains $S_j$. This implies
there is a subscheme of $H_0 \cap X_0$ of length at least
$2g-2-l+s-1$ having support contained in the finite set $(\Lambda _0
\cap X_0)\cup \{ S_2, \cdots, S_s \}$. Consider the normalization
$n: X' \rightarrow X_0$ then $n^{-1}((\Lambda _0 \cap X_0)\cup \{
S_2, \cdots, S_s \})$ is a finite subset $Z$ of $X'$. There are
finitely many effective divisors of degree at most $2g-2$ having
support in $Z$. The embedding $X_0 \subset \mathbb{P}^{g-1}$ gives
rise to a morphism $m:X' \rightarrow \mathbb{P}^{g-1}$ and since the
image has degree $2g-2$ this morphism corresponds to a
subvectorspace $V$ of dimension $g$ of the space of global sections
of an invertible sheaf $\mathcal{L}$ of degree $2g-2$ on $X'$. For
each $P\in X'(\mathbb{R})$ we obtain a section in that vectorspace
defining a divisor containing $P$ and having a subdivisor of degree
at least $2g-3-l+s$ having support contained in $Z$. Since
$X'(\mathbb{R})$ has infinitely many points there exists an
effective divisor $D_0$ of degree at least $2g-3-l+s$ having support
on $Z$ satisfying the following property. There exist infinitely
many points $P$ on $X'(\mathbb{R})$ such that there exist a section
$s$ of $\mathcal{L}$ whose associated divisor $D$ contains $D_0+P$.
For those divisors $D-D_0$ is an effective divisor corresponding to
a global section of the invertible sheaf $\mathcal{L}(-D_0)$ on $X'$
of degree at most $1+l-s<1+k-s=k'$. Since there exist infinitely
many points $P$ on $X'(\mathbb{R})$ belonging to the divisor of a
global section of $\mathcal{L}(-D_0)$ it implies
$h^0(\mathcal{L}(-D_0))\geq 1$. This contradicts the assumptions on
$X'$.

\end{document}